%% file: agt-1-26.tex
\documentclass{gtart}
\input agtout

\lognumber{26}
\volumenumber{1}
\volumeyear{2001}
\papernumber{26}
\pagenumbers{503}{518}
\received{6 August 2001}
\revised{20 September 2001}
\accepted{21 September 2001}
\published{22 September 2001}

\usepackage{amsmath,amssymb}
\usepackage{psfrag,graphicx}
\usepackage{amsfonts}

\newtheorem{thm}{Theorem}[section]
\newtheorem{lemma}[thm]{Lemma}

\newtheorem{prop}[thm]{Proposition}

\theoremstyle{definition}
\newtheorem{defn}[thm]{Definition}
\newtheorem{example}[thm]{Example}

\newcommand{\bC}{{\mathbb {C}}}

\newcommand{\bR}{{\mathbb {R}}}
\newcommand{\bZ}{{\mathbb {Z}}}
\newcommand{\cC}{{\mathcal C}}
\newcommand{\cL}{{\mathcal L}}
\newcommand{\cS}{{\mathcal S}}
\newcommand{\cR}{{\mathcal R}}
\newcommand{\cZ}{{\mathcal Z}}

\newcommand{\bdy}{\partial}
\newcommand{\del}{\partial}
\newcommand{\wtilde}{\widetilde}

\newcommand{\vphi}{\varphi}

\begin{document}
\title{Generalized symplectic rational blowdowns}
\author{Margaret Symington}
\address{School of Mathematics\\Georgia Institute of 
Technology\\Atlanta, GA 30332, USA}
\email{msyming@math.gatech.edu}

\begin{abstract}
We prove that the generalized rational blowdown, 
a surgery on smooth $4$-manifolds,
can be performed in the symplectic category. 
\end{abstract}
\asciiabstract{We prove that the generalized rational blowdown, a
surgery on smooth 4-manifolds, can be performed in the symplectic
category.}

\primaryclass{57R17} 
\secondaryclass{57R15, 57M50}
\keywords{Symplectic surgery, blowdown}

\maketitle

\section{Introduction}

Surgery techniques 
are essential tools for understanding the topology of manifolds.
For smooth manifolds the rational blowdown surgery, 
introduced by Fintushel and 
Stern, is particularly useful because one can calculate how the
Donaldson and Seiberg-Witten invariants change when the surgery is 
performed~\cite{FintStern.blowdown}.
For instance,
Fintushel and Stern~\cite{FintStern.blowdown} used it to calculate the 
Donaldson and Seiberg-Witten invariants
of simply connected elliptic surfaces and to
construct an interesting family of simply connected smooth $4$-manifolds
$Y(n)$ not homotopy equivalent to any complex surface.
This surgery can also be performed in the symplectic 
category~\cite{Sym.blowdowns}, and thereby helps demonstrate the vastness of
the set of symplectic $4$-manifolds.
In particular, the aforementioned
$Y(n)$, as well as an infinite family of exotic K3 
surfaces~\cite{GompfMrowka} 
($4$-manifolds that are homeomorphic but not diffeomorphic to 
a degree $4$ complex hypersurface in $\bC P^3$), all admit symplectic
structures~\cite{Sym.blowdowns}.

The rational blowdown surgery 
amounts to removing a neighborhood of a linear chain
of embedded spheres whose boundary is the lens
space $L(n^2,n-1)$, $n\ge 2$
and replacing it with a
{\it rational ball} 
(manifold with the same rational homology as a ball), also with boundary
$L(n^2,n-1)$.
This has the effect of reducing the dimension of the second homology of
$M$ at the expense of possibly 
complicating the fundamental group.
The surgery gets its name from the well-known process of blowing down
a $-1$ sphere (the case $n=1$) in which one replaces a
tubular neighborhood of a sphere of self-intersection $-1$ by a
$4$-ball.

In fact, there are other lens spaces that bound rational balls:
$L(n^2,nm-1)$, $n, m\ge 1$ and relatively prime~\cite{CassonHarer}.
Therefore one can define a broader class of rational blowdowns, so called
{\it generalized rational blowdowns}.
Park~\cite{Park.gblowdowns} extended Fintushel and Stern's calculations,
showing how a generalized rational
blowdown affects the Donaldson and Seiberg-Witten invariants of a smooth
$4$-manifold.
Here we show that even the
generalized rational blowdown can be performed
in the symplectic category.

Specifically,
given any pair of relatively prime integers $y,x$, $x\ne 0$, the fraction
$\frac{y}{x}$ has a {\it negative continued fraction expansion}
\begin{equation}
 b_1 - \left(\frac{1}{b_2 - \frac{1}{\cdots - \frac{1}{ b_k} }}\right) = \frac{y}{x}
\end{equation}
which is unique if one assumes that $b_j\ge 2$ for all $j\ge2$.
The shorthand for this continued fraction expansion is
$[b_1,b_2,\ldots, b_k]$.

\begin{figure}
\small
\begin{center}
\psfrag{b1}[][]{ $b_1$}
\psfrag{b2}[][]{ $b_2$}
\psfrag{b3}[][]{ $b_3$}
\psfrag{bn1}[][]{ $b_{k-1}$}
\psfrag{bn}[][]{ $b_k$}
\includegraphics[angle=0,width=3.0in]{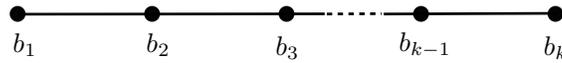}
\caption{Plumbing diagram for $C_{n,m}.$}
\label{plumb.fig}
\end{center}
\end{figure}

\begin{defn}
For any relatively prime $n\ge 2, m\ge 1$, 
let $C_{n,m}$ be a closed tubular neighborhood of the union of spheres 
$\{S_j\}_{j=1,\ldots k}$ in the plumbing of disk bundles represented by
the diagram in Figure~\ref{plumb.fig},
where the $b_j$ satisfy 
$[b_1,b_2,\ldots, b_k]=\frac{n^2}{nm-1}$ and $b_j\ge 2$ for all $j$.
\end{defn}

The spheres in $C_{n,m}$ have
the following intersection pattern:
\begin{equation}
\left\{
\begin{array}{lll}
S_j\cdot S_{j+1}&=1  &   {\rm for}\ \  j=1,\ldots k-1, \\
S_j\cdot S_j&=-b_j & {\rm and} \\
S_i\cdot S_j&=0 & {\rm otherwise}. 
\end{array}
\right.
\end{equation}

The fact that $S_j\cdot S_j=-b_j\le -2$ for each $j$ implies
that the intersection form of $C_{n,m}$ is negative definite.
The boundary of $C_{n,m}$ is the lens space $L(n^2,nm-1)$
which bounds a rational homology
ball $B_{n,m}$~\cite{CassonHarer,Park.gblowdowns}.

\begin{defn} If there is an embedding $\psi\co C_{n,m}\rightarrow M$, 
then the generalized rational
blowdown of $M$ along the spheres $\psi(\cup_{i=1}^k S_i)$ is 
\begin{equation}
\wtilde M :=(M- \psi(\cup_{i=1}^k S_i))\cup_\phi B_{n,m}
\end{equation}
 where $\phi$ is
an orientation preserving diffeomorphism of a collar
neighborhood of the boundary $L(n^2,nm-1)$.
\end{defn}

\begin{thm} \label{gbd.thm} Suppose 
$\wtilde M=(M- \psi(\cup_{i=1}^k S_i))\cup_\phi B_{n,m}$ 
is the generalized
rational blowdown of a smooth $4$-manifold $M$ along spheres 
$\psi(\cup_{i=1}^k S_i)$.
If $M$ admits a symplectic structure
for which the spheres are symplectic,
then the diffeomorphism $\phi$ can be chosen so that
$\wtilde M$ admits a symplectic structure 
induced from the symplectic structures on $M$ and $B_{n,m}$.
\end{thm}

The essence of the proof, as for the case $m=1$, is in 
our choice of symplectic models for the spaces 
$C_{n,m}$ and  $B_{n,m}$.
By a version of the symplectic neighborhood theorem any neighborhood of
symplectic spheres that is diffeomorphic to $C_{n,m}$ has a 
neighborhood symplectomorphic to a toric model space.
(A symplectic manifold is {\it toric} if it is equipped with an effective
Hamiltonian $T^n$ action.)

The new ingredient in this paper is a set of 
symplectic representatives for the
rational balls $B_{n,m}$ for all $m\ge 1$.
These representatives  are 
toric near the boundary and can be chosen to \lq\lq fit\rq\rq \ 
a collar neighborhood of the boundary of $C_{n,m}$.
We present the $B_{n,m}$ as the total space of 
a singular Lagrangian fibration with two types of singular fibers:
a one parameter family of circle fibers and one isolated nodal fiber --
a sphere with one positive self-intersection.
In the language of Hamiltonian integrable systems the singularity of
the nodal fiber is a focus-focus singularity.
Our nodal fiber is the Lagrangian analog of the singular fibers that appear
in Lefschetz fibrations of symplectic $4$-manifolds.

{\bf Acknowledgments}\qua The author thanks Eugene Lerman for suggesting
helpful references, in particular the work of Nguyen Tien Zung, and 
thanks Nguyen Tien Zung in turn for mentioning Vu Ngoc San's work.
The author also thanks an anonymous referee who pointed out a minor error.

The author is grateful for the support of an NSF post-doctoral 
fellowship, DMS9627749.

\section{Background}

Our objective is to control the symplectic structure of collar
neighborhoods of the boundaries of
the spaces involved in our surgery: $C_{n,m}$ and
$B_{n,m}$.
We do this by presenting them as the total spaces of singular Lagrangian
fibrations.
The space $C_{n,m}$ itself and a collar neighborhood of the boundary of 
$B_{n,m}$ admit singular Lagrangian fibrations equivalent to the fibration
defined by the moment
map for a Hamiltonian torus action.
An important feature of these fibrations is that, at least near the
boundary, the base classifies the neighborhood up to 
fiberwise symplectomorphism 
(cf.~\cite{BouMol.fibr,Zung.II}).

\begin{defn} A Lagrangian fibration $ \pi\co (M^{2n},\omega)\rightarrow B^n$
is a locally trivial fibration such that 
$\omega|_{\pi^{-1}(b)}= 0$ for each $b\in B$ (i.e. such that each fiber is
a Lagrangian submanifold).
\end{defn}

The Arnold-Liouville theorem guarantees that if the fibers of a Lagrangian
fibration are closed (compact, without boundary) and connected
then they must be $n$-tori with neighborhoods equipped with canonical
coordinates: action-angle coordinates.
The local action coordinates supply the base $B$ with an integral 
affine structure, i.e. an atlas 
$\Phi_j\co U_j\rightarrow \bR^n$ 
with the maps $\Phi_j\Phi_k^{-1}|_{U_j\cap U_k}
\in GL(n,\bZ)$.

It is easy to see that in dimension $4$ (with $n=2$) the fibers must be tori:
the Lagrangian condition
implies the existence of an isomorphism between the normal and tangent bundles
defined via the symplectic form; then, since the normal bundle of a fiber 
must be trivial, we have that Euler characteristic of the tangent bundle
is $0$.  Because the fiber of a locally trivial fibration of an oriented
manifold is orientable, the fiber must be a torus.

We now expand our definition of a Lagrangian fibration to include singular
fibers: one parameter families of circle fibers, isolated points and
isolated nodal fibers (spheres with one positive transverse intersection).
These singular fibrations
are examples of Lagrangian fibrations with topologically stable and 
non-degenerate singularities such as 
arise in integrable systems~\cite{Zung.II}. 
In the spirit of holomorphic fibrations and smooth Lefschetz fibrations,
and for simplicity of exposition, we often suppress the word singular.
We assume throughout that fibers are connected and that the generic
fibers are closed manifolds.

Near the circle and point fibers the fibration is equivalent to one
coming from the moment map for a torus action.
Therefore, the integral affine structure on the image of the 
regular fibers, $B_0\subset B$, extends to each connected component 
of the image of the circle fibers.
These components meet at the vertices of $B$, the images of the
point fibers.
The images of nodal fibers are isolated interior points of $B$.

To understand the base $B$ in each of our examples, we view $B$ (or
part of it) as a subset of $\bR^2$. 
We always assume that $\bR^2$ is equipped with the integral
affine structure coming from the standard lattice generated by the
vectors $(1,0)$ and $(0,1)$.
It is important to note that there are two different classes of lines in
this integral affine space: rational and irrational (as determined by
the slope of the line).
Indeed, a vector $v$ directed along a line of rational
slope has an affine length $\alpha\in\bR^+$ 
defined by $v=\alpha u$ for a primitive
integral vector $u$, while a vector directed along a line of irrational
slope does not have a well-defined length.
By an {\it integral polygon} in $\bR^2$ we mean one whose edges define
vectors of the form $\alpha u$ where $\alpha\in \bR^+$ and $u$ is a
primitive integral vector, or alternatively, one whose edges all have
well-defined affine lengths.

We now review a few facts that facilitate reading the topology 
of a Lagrangian fibered symplectic $4$-manifold 
$\pi\co (M,\omega)\rightarrow U$
from the base $U$ when $U$ coincides with a moment map image.
The reader interested in more detail on the topology of a toric symplectic
manifold should consult~\cite{Audin.torus}.
Throughout our discussion, {\it neighborhood} 
refers to a tubular neighborhood, $(p_1,p_2)$ are Euclidean coordinates on 
$\bR^2$ and $(q_1,q_2)$ are circular coordinates on $T^2$.

\begin{enumerate}

\item
A simply connected 
open domain $U\subset \bR^2$ defines the open symplectic manifold
$(U\times T^2, dp\wedge dq)$.

\item An open 
neighborhood $U$ of a point in the boundary of a closed half-plane
in $\bR^2$
defines the smooth manifold  $S^1\times D^3$
that is symplectomorphic to a neighborhood of $\{(z_1,z_2)\,|\,
0<|z_1|<\alpha,|z_2|=0\}\subset (\bC^2,\frac{1}{2i}d{\overline z}\wedge dz)$
for some $\alpha\in \bR^+$.
The symplectic structure $dp\wedge dq$ defined on the preimage of
${\rm int} \ U$ extends to the circles that live over the points in 
$\bdy U\cap U$.
If the half space is  
bounded by the line $\{(p_1,p_2)\,|\,p_2=\frac{m}{n} p_1\}\subset\bR^2$ 
then the circle fibers are quotients
 $T^2/(q_1,q_2)\sim(q_1-mt,q_2+nt)$, $t\in \bR$
of the tori living
over points in  $U$ with $p_2=\frac{m}{n} p_1$.

\item
A neighborhood $U$ of a vertex in a convex integral polygon
defines a symplectic $4$-ball if and only if
the primitive integral vectors $u,v$ that define the directions of
the adjacent edges satisfy
$|u\times v|=1$.  (Here $\times$ is the cross product in $R^3$ restricted to
$R^2$, thus yielding a scaler.)
The preimage of the vertex is a point.
If $|u\times v|=n\ge 2$ then $U$ defines
a neighborhood of an orbifold singularity.

\item
A neighborhood $U$ 
of an edge $E$ in a convex polygon defines a neighborhood of a sphere.
Specifically, consider an edge $\alpha w$, with $\alpha\in\bR^+$ and $w$ a
primitive integral vector.
Suppose
$u,v$ are the primitive integral vectors based at the endpoints of $E$ 
that define (up to scaling) the left and right adjacent edges.
Then the sphere has area $\alpha$ and self-intersection
$u\times v$.  See Figure~\ref{sphere.fig}.

\begin{figure}
\small
\begin{center}
\psfrag{u}[][]{ $u$}
\psfrag{v}[][]{ $v$}
\psfrag{w}[][]{ $w$}
\includegraphics[angle=0,width=3.0in]{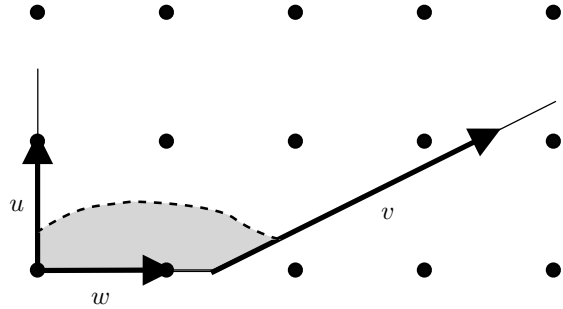}
\caption{Neighborhood of a sphere of self-intersection $-2$ and area
$\frac{3}{2}$.}
\label{sphere.fig}
\end{center}
\end{figure}

\item 
If $U\subset \bR^2$ defines a toric symplectic manifold, then 
for any $A\in GL(2,\bZ)$ and $b\in\bR^2$, $A(U)+b$ defines the same
symplectic manifold (with a different torus action if $A$ is not the
identity). 

\end{enumerate}

\section{Symplectic models} \label{models.sec}

In this section we provide symplectic models for the cone on a lens space,
neighborhoods of certain linear chains of spheres, the neighborhood
of a nodal fiber, and rational balls.  
We give the descriptions in terms of diagrams in $\bR^2$ that correspond
to images of moment maps when a global torus action can be defined.

The examples we present here are the building blocks for our constructions
and are essential for the proof of Theorem~\ref{gbd.thm}.

\subsection{Toric models}

\begin{example} \rm {Cone on a lens space $L(n,m)$.} \label{cone.ex}

Consider the following subset of $\bR^2$:
\begin{equation}V_{n,m}=\{p_1\ge 0\}\cap\{p_2\ge \frac{m}{n}p_1\}\cap\{p_2>0\}
\end{equation}
and the (singular) Lagrangian fibered symplectic manifold
$\pi\co (M,\omega)\rightarrow V_{n,m}$ it defines.
Figure~\ref{cone.fig} shows $V_{n^2,nm-1}$, the case we are interested in.

To see that $M$ is a cone on a lens space,
recall that $L(n,m)$ can be decomposed as 
the union of two solid tori glued together via a map $\phi$ 
of their boundaries such that $\phi_*\mu_2= -m\mu_1 + n\lambda_1$ 
where $\mu_i,\lambda_i$ are meridinal and longitudinal cycles 
on the solid torus boundaries.

\begin{figure}
\small
\begin{center}
\psfrag{nm}[][]{ $nm-1$}
\psfrag{n2}[][]{ $n^2$}
\psfrag{Vn2nm}[][]{ $V_{n^2,nm-1}$}
\includegraphics[angle=0,width=3.0in]{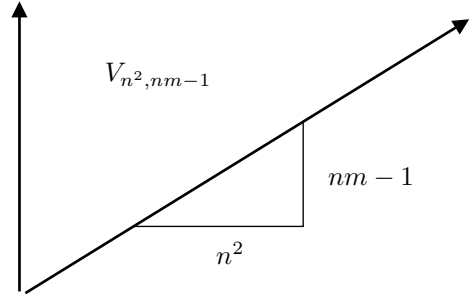}
\caption{Cone on the lens space $L(n^2,nm-1)$.}
\label{cone.fig}
\end{center}
\end{figure}

For any $t>0$, consider the
$3$-manifold in $M$ that is the preimage of $\{p_2=t\}\cap V_{n,m}$;
decompose it as the union of preimages $P_1\cup P_2$
where $P_1$ is the preimage of $\{p_1\le ct,\ p_2=t\}\cap V_{n,m}$
and $P_2$ is the preimage of $\{p_1\ge ct,\ p_2=t\}\cap V_{n,m}$
for some $0<c<\frac{n}{m}$.  
Then $P_1, P_2$ are a solid tori with meridians whose tangent vectors
are $\frac{\del}{\del q_1}$ and
$-m\frac{\del}{\del q_1}+n\frac{\del}{\del q_2}$ respectively, thereby showing
the $3$-manifold is 
$L(n,m)$.  
Letting $t$ vary we get $L(n,m)\times (0,\infty)$.
\end{example}

There was nothing special about the choice of 
$\{p_2=t\}\cap V_{n,m}$ to define the lens space; we could have
used any arc smoothly embedded in $V_{n,m}$ with one endpoint on each of
the edges of $V_{n,m}$.
However, by choosing an arc $\gamma$ transverse to the vector field
$p_1 \frac{\del}{\del p_1}+p_2 \frac{\del}{\del p_2}$
we get an induced contact structure (completely non-integrable
$2$-plane field, cf.~\cite{Etnyre.convexity}) on the lens space.
This contact structure is defined as the kernel of the $1$-form
$\iota_X\omega|_{\pi^{-1}(\gamma)}$ 
where $X$ is the unique vector field on $M$ which is given
by $p_1\frac{\del}{\del p_1} + p_2\frac{\del}{\del p_2}$ in the local
coordinates $(p,q)$ on $\pi^{-1}({\rm int\ } V_{n,m})$.
The contact structure is independent of the choice of the transverse arc
$\gamma$.

\begin{example}\rm {Negative definite chains of spheres.}\label{nbhd.ex}

Here we define a neighborhood of a chain of spheres by a neighborhood
of the piecewise linear boundary of a domain in $\bR^2$.
See Figure~\ref{nbhd.fig} for an example.

Let $\{x_j\}_{j=0}^{k}$ be a set of points in $\bR^2$ and
$\{u_j\}_{j=0}^{k+1} $ a set of primitive integral vectors such that

\begin{itemize}
\item $\alpha_j u_j=x_{j}-x_{j-1}$ with $\alpha_j\in\bR^+$
for each $1\le j\le k$,
\item $u_{j}\times u_{j+1}=1$ for each $0\le j\le k$, and
\item $u_{j+1}\times u_{j-1}=S_j\cdot S_j$ for each $1\le j\le k$.
\end{itemize}
Let $X$ be the convex hull of the points $\{x_j\}_{j=0}^{k}$ and all
points $x$ such that $x_0-x=\alpha u_0$ or $x-x_{k}=\alpha u_{k+1}$
for some $\alpha>0$.

\begin{figure}
\small
\begin{center}
\psfrag{u0}[][]{ $u_0$}
\psfrag{u1}[][]{ $u_1$}
\psfrag{u2}[][]{ $u_2$}
\psfrag{u3}[][]{ $u_3$}
\psfrag{u4}[][]{ $u_4$}
\psfrag{x0}[][]{ $x_0$}
\psfrag{x1}[][]{ $x_1$}
\psfrag{x2}[][]{ $x_2$}
\psfrag{x3}[][]{ $x_3$}
\psfrag{U}[][]{ $W$}
\psfrag{X}[][]{ $X$}
\includegraphics[angle=0,width=3.0in]{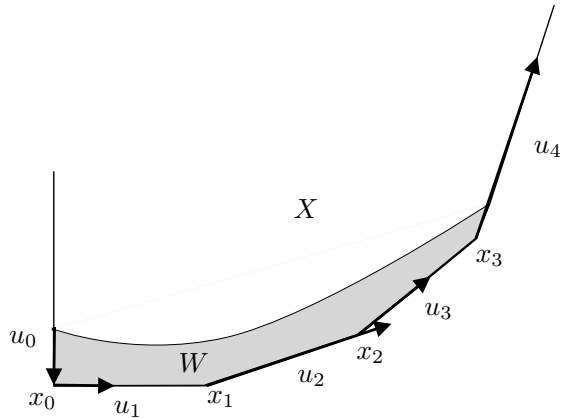}
\caption{Neighborhood of spheres.}
\label{nbhd.fig}
\end{center}
\end{figure}

Then $X$ defines a Lagrangian fibered symplectic manifold 
$(M,\omega)\rightarrow X$ such that each finite
edge, defined by the vector $x_{j}-x_{j-1}$ for some $j$, is the image
of a sphere $S_j$.
The area of each sphere $S_j\subset M$ is $\alpha_j$ and
for each $1\le j\le k-1$, 
$S_j$ intersects $S_{j+1}$ once positively and transversely.
The convexity of $V$ corresponds to the negative definiteness of the 
intersection form of $M$.

Let $W$ be any closed neighborhood in $X$ of the finite edges defined by the
$x_j$.
Then $W$ defines a singular Lagrangian fibration of 
a closed toric neighborhood of 
spheres $S_1, \ldots S_k$ in $M$.
(We interpret the points in 
$\bdy W \cap {\rm int}\  X$ as the images of tori, not circles.)

\end{example}

A variation of the symplectic neighborhood theorem states that
the germ of the neighborhood of a linear chain of spheres
is determined up to
symplectomorphism by the areas of the spheres 
and the intersection form.
(An explanation of how Moser's method would be applied in
this case is provided in~\cite{McR.darboux}.)
Therefore, given any symplectic manifold $(M,\omega)$ containing
a smoothly embedded copy of $C_{n,m}$ 
as a neighborhood of symplectic spheres, 
we can choose an $X=X_{n,m}$ and a $W_{n,m}\subset X_{n,m}$
small enough that $W_{n,m}$  
defines a Lagrangian fibered symplectic manifold that symplectically
embeds in and is diffeomorphic to the embedded copy of $C_{n,m}$.
Therefore, we simply assume that $C_{n,m}$ is symplectically
embedded in $M$ and Lagrangian
fibers over $W_{n,m}$.
We also assume, without loss of generality, that 
the boundary of $C_{n,m}$ has an induced
contact structure equivalent to the one described in Example~\ref{cone.ex}
when the lens space is $L(n^2,nm-1)$.

\subsection{Neighborhood of a nodal fiber}\label{nodal.sec} \

Nodal fibers appear as singular fibers in numerous integrable systems 
including the spherical pendulum (cf.~\cite{Duist.actionangle, Zung.focus}).
As noted by Zung~\cite{Zung.focus}, a
simple model for a Lagrangian fibered neighborhood of 
a nodal fiber is a self-plumbing of the zero section of 
$(T^*S^2,\omega={\rm Re}\  dz_1\wedge
dz_2)$.
Indeed, glue a neighborhood of $(0,0)\subset \bC^2$ to a neighborhood
of $(\infty,0)$ by the symplectomorphism
$(z_1,z_2)\rightarrow (z_2^{-1},z_1z_2^2)$.
Projecting to $\bR^2=\bC$ by the map $z_1z_2$ gives the desired Lagrangian
fibration.

\begin{lemma} \label{germ.lem} The germ of a symplectic neighborhood of 
a Lagrangian nodal fiber is unique up to symplectomorphism.
\end{lemma}

\begin{proof}
The lemma follows from the Lagrangian neighborhood theorem by pulling the
symplectic structure of a nodal fiber neighborhood back to a neighborhood
of the zero section of $T^*S^2$ via an immersion.
\end{proof}

Let $\pi\co (N,\omega)\rightarrow B$ 
be a Lagrangian fibered neighborhood of a nodal fiber 
with $B$ a disk and $b_0\in B$ the image of the nodal fiber.
The Arnold-Liouville theorem implies that 
$B-b_0$ is equipped with an integral affine structure.
In particular, $T(B-b_0)$ has a flat connection.
The topological monodromy 
\begin{equation}
A=\left(  
\begin{array}{cc} 
1 & 1 \\  
0 & 1 
\end{array} 
\right)
\end{equation}
of the torus fibration over $B-b_0$ and the Lagrangian
structure of the fibration forces the same monodromy in the induced
flat connection on $T(B-b_0)$.
Therefore no embedding of 
$B$ into $\bR^2$ preserves the (integral) affine structure.
However, we can find a map that is an isomorphism almost everywhere.

Indeed, $B-b_0$ must be  
isomorphic to a neighborhood of the puncture  in a punctured plane
with integral affine structure and monodromy
$A.$
Specifically, let $X$ be the universal cover of $\bR^2-0$ with the
affine structure lifted from $\bR^2$ and polar coordinates $(r,\theta)$,
$-\infty<\theta<\infty$.
With $p=(p_1,p_2)$ the Euclidean coordinates on $\bR^2$, 
we can also identify points
in $X$ by $(p,n)$ where $n=\left[\frac{\theta}{2\pi}\right]$.
Let $V_n\subset X$, $n\ge 1$ 
be defined by $0 < \theta < 2n\pi + \frac{\pi}{2}$.
Define the sectors  $S_n, S_0\subset V_n$ by 
$2n\pi < \theta < 2n\pi + \frac{\pi}{2}$ 
and $0 < \theta < \frac{\pi}{4}$  respectively.
Now glue the sector $S_n$ to the sector $S_0$ 
via the map that, with respect to the
labeling $(p,n)$, sends the point $(p,n)$ to $(Ap,0)$.
Call the resulting manifold $P_n$.

\begin{lemma}\label{plane.lem}
Each $P_n$ defines a Lagrangian fibration 
$\pi \co (M_n,\omega_n)\rightarrow P_n$
that is unique up to fiberwise symplectomorphism.
\end{lemma}

\begin{proof}
We can construct a Lagrangian torus fibration with base $P_n$ 
as follows:
equip  $V_n\times T^2$ with coordinates $(p,q,n)$ and symplectic form 
$dp\wedge dq$ where $q=(q_1,q_2)$ are coordinates on the torus.
Now glue $S_n\times T^2$ to $S_0\times T^2$ via the symplectomorphism that
sends $(p,q,n)$ to $(Ap,A^{-T}q,0)$. 
The resulting manifold is $M_n$; forgetting the torus
coordinates $q$ gives the desired fibration over $P_n$.
This Lagrangian fibration is uniquely defined by the base
because $P_n$ has the homotopy type
of a $1$-dimensional manifold (\cite{Duist.actionangle}).
\end{proof}

This lemma is clearly still true 
if we replace $P_n$ with a neighborhood
$U_n$ of the puncture in $P_n$.
Furthermore, 
two such neighborhoods $U_n$, $U'_n$ define a symplectically equivalent
Lagrangian fibrations if and only if they are integral affine isomorphic. 
Note that in terms of the coordinates used in the proof of 
Lemma~\ref{plane.lem} the vector field $\frac{\del}{\del q_2}$ on
$V_n\times T^2$ descends to a well defined vector field on $M_n$ which
for simplicity we also call  $\frac{\del}{\del q_2}$.

\begin{lemma}
Let $\pi\co (N,\omega)\rightarrow B$ be a singular Lagrangian fibration 
with one singular fiber, a nodal fiber with image $b_0\in B$ where $B$ is
a disk.
The punctured disk $B-b_0$ is affine isomorphic to a neighborhood of the
puncture in $P_1$
and $N-\pi^{-1}(b_0)$ symplectically embeds in $M_1$ as the preimage of
some $U_1$.
The vanishing cycle of the nodal fiber is
the cycle represented by an integral curve of the vector field 
$\frac{\del}{\del q_2}$ on $M_1$.
\end{lemma}

\begin{proof}
One can see that $n=1$ in one of two ways: 
Duistermaat~\cite{Duist.actionangle} calculated
explicit action coordinates in a neighborhood of a nodal
fiber -- on the complement of the fibration over a ray based
at $b_0$.
In other words, he found the aforementioned isomorphism.
Alternatively, if the boundary  of
$B$ is chosen to be transverse to rays emanating from $b_0$ then
the boundary of $N$
is equipped with a contact structure induced from the symplectic
structure on $N$.
Because this contact structure is fillable, it must be tight 
(cf.~\cite{Etnyre.convexity}), but this can happen only if $n=1$; otherwise the
structure would be overtwisted~\cite{Lerman.ccuts}.

The vanishing cycle is in the class of the eigenvector of the monodromy matrix
for the torus bundle fibering over $B-b_0$.  
Appealing to the model
$M_1$ constructed in the proof of Lemma~\ref{plane.lem}, we see  
this is the eigenvector
of $A^{-T}$, namely $\frac{\del}{\del q_2}$.
\end{proof}

In $P_1$, with coordinates chosen as above, we call 
the line in the base defined by the vector $(1,0)$ the {\it eigenline}.
It is the only well defined line that passes through the puncture.

Two neighborhoods $(N_0,\omega_0)$, $(N_1,\omega_1)$ 
of nodal fibers that are Lagrangian
fibrations over the same base $B$ need not be 
fiberwise symplectomorphic.
Indeed, there is a Taylor series invariant of the Lagrangian fibration -- an
element of $\bR[[X,Y]]_0$, 
the algebra of formal power series in two variables with vanishing
constant term -- that classifies such a neighborhood up to {\it fiberwise}
symplectomorphism~\cite{San.focus}.
However, we are only interested in classifying the neighborhood up
symplectomorphism.

\begin{lemma} \label{nbhd.lem}
Two neighborhoods $(N_0,\omega_0)$, $(N_1,\omega_1)$ 
of nodal fibers that are singular Lagrangian
fibrations over the same base $B$ are symplectomorphic.
\end{lemma}

\begin{proof}
Let $\cS_0,\cS_1\in \bR[[X,Y]]_0$ be the Taylor series that classify
the germs of the neighborhoods of the singular fibers in $N_0,N_1$.
Following  San~\cite{San.focus}, we can use two functions 
$S_0,S_1\in \cC^\infty(\bR^2)$ whose Taylor series are $\cS_0,\cS_1$
to construct model Lagrangian fibered symplectic neighborhoods 
equivalent to $N_0,N_1$.
We can choose $S_0,S_1$ to be equal outside 
of a small neighborhood $V$ of the origin and
then choose a smooth family of functions $S_t$ that vanish at the identity,
connect $S_0$ and $S_1$, and are equal to $S_0$ and $S_1$
outside of $V$.
Using these functions 
we can construct a $1$-parameter family of Lagrangian fibered
neighborhoods $(N_t,\omega_t)$.
It is then easy to define a $1$-parameter family of
diffeomorphisms $\vphi_t\co N_0\rightarrow N_t$ such that $\vphi_0$ is the 
identity and $\vphi_t^*\omega_t=\omega_0$ on the complement of a smaller
neighborhood of the nodal fiber.
Because the induced symplectic forms $\vphi_t^*\omega_t$ 
are all cohomologous a Moser argument completes the
proof.
\end{proof}

\subsection{Symplectic rational balls}
\label{rball.sec} 

To prove Theorem~\ref{gbd.thm} we need symplectic models for the
rational balls $B_{n,m}$ whose boundaries are the lens spaces $L(n^2,nm-1)$.
We do this by defining Lagrangian fibrations
$\pi\co (B_{n,m},\omega_{n,m})\rightarrow U_{n,m}$
with two types of singular fibers: a one parameter family of
circle fibers and one nodal fiber.

First note that in our construction of a model neighborhood of a nodal
fiber we can make a different choice of coordinates, with
respect to which the eigenline is in the $(n,m)$ direction in $\bR^2$
and the vanishing cycle is in the class of an integral curve of 
$-m\frac{\del}{\del q_1}+n\frac{\del}{\del q_2}$.
(Here $m$ and $n$ are relatively prime integers.)

Now let $A_{n,m}$ be a space diffeomorphic to a closed half-plane in $\bR^2$
and such that:

\begin{itemize}
\item there is a distinguished point $p_t\in {\rm int}\  A_{n,m}$ such
that $A_{n,m}-p_t$ is equipped with an affine integral
structure and the monodromy around $p_t$ is 
$\left(  
\begin{array}{cc} 
1& 1 \\  
0 & 1
\end{array} 
\right)$;
\item the eigenline through $p_t$ intersects the boundary in a point $p_0$; 
and 
\item $A_{n,m}$ minus the line segment $L_t$
connecting $p_0$ and $p_t$ is affine
isomorphic to the following domain in $\bR^2$:
\begin{equation}
\{(p_1,p_2)\,|\, p_1\ge 0, p_2 \ge \frac{nm-1}{n^2}, p_2>0 \}
\end{equation}
minus the line segment connecting the points $(0,0)$ and $(tn,tm)$ for
some $t>0$.

\end{itemize}

\begin{figure}
\small
\begin{center}
\psfrag{nm}[][]{ $nm-1$}
\psfrag{n2}[][]{ $n^2$}
\psfrag{n}[][]{ $n$}
\psfrag{m}[][]{ $m$}
\psfrag{pt}[][]{ $p_t$}
\psfrag{Unm}[][]{ $U_{n,m}$}
\psfrag{Anm}[][]{ $A_{n,m}$}
\includegraphics[angle=0,width=3.0in]{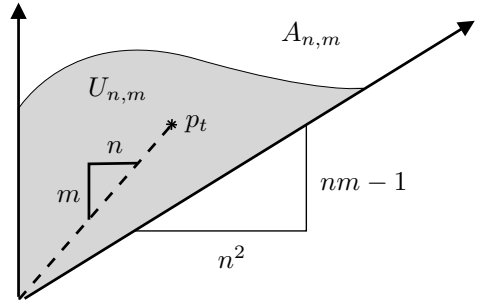}
\caption{Rational ball with boundary $L(n^2,nm-1)$.}
\label{rball.fig}
\end{center}
\end{figure}

Let $U_{n,m}\subset A_{n,m}$ be a closed neighborhood of $L_t$ (which
necessarily contains a connected segment of $\bdy A_{n,m}$).  
In Figure~\ref{rball.fig} we show the image of $U_{n,m}-L_t\subset A_{n,m}-L_t$
in $\bR^2$ under the aforementioned isomorphism. 

\begin{lemma}
$U_{n,m}$ is the base of a (singular) Lagrangian fibration of 
the rational ball $B_{n,m}$. 
\end{lemma}

\rk{Remark}
In this description we understand that the preimage of points in 
$\bdy U_{n,m}\cap  {\rm int}\  A_{n,m}$ are tori so that $U_{n,m}$ 
defines a manifold with boundary.
The image of the boundary is the closure of 
$\bdy U_{n,m}\cap  {\rm int}\  A_{n,m}$ in $A_{n,m}$.

\begin{proof} 
Because it is homotopic to a $1$-manifold, 
$U_{n,m}-p_t$ defines a unique
Lagrangian fibration $\pi\co M_0\rightarrow U_{n,m}-p_t$ 
with $\pi^{-1}(b)$ a circle for each 
$b\in \bdy U_{n,m}\cap \bdy A_{n,m}$ (cf.~\cite{BouMol.fibr,Zung.II}).

An open 
neighborhood of $p_t\subset U_{n,m}$ 
is the base of a singular Lagrangian
fibration of a neighborhood of a nodal fiber as in Section~\ref{nodal.sec}.
Therefore we can glue a neighborhood of a nodal fiber into $M_0$ with
a fiber-preserving symplectomorphism to get a symplectic manifold $M$
fibering over $U_{n,m}$. 

To see that $M$ is a rational ball it
suffices to note that it is homotopy equivalent to
the preimage of an embedded arc connecting the boundary of 
$U_{n,m}$ and $p_t$.
This preimage is homeomorphic to the 
the space obtained from
$T^2\times[0,1]$ by collapsing all $(1,0)$ curves on  $T^2\times \{0\}$ (to
get the circle fiber over the boundary point)
and a $(-m,n)$ curve on $T^2\times \{1\}$ (to get the nodal fiber).
Because $n\ne 0$, we see $H_1(M,\bR)=H_2(M,\bR)=0$ and $\pi_1(M)=\bZ_n$.

Finally, $M=B_{n,m}$ because its boundary is the lens space $L(n^2, nm-1)$
as can be seen by comparing Figures~\ref{cone.fig} and~\ref{rball.fig}:
a collar neighborhood of the boundary of $M$ projects to 
a subset of $U_{n,m}$ which is clearly isomorphic to a one sided neighborhood
of an arc connecting the two boundary components of $V_{n^2,nm-1}$. 
(See Example~\ref{cone.ex}.)
\end{proof}
 
\begin{prop} \label{ball.prop}
For a given $U_{n,m}$, 
the rational ball $B_{n,m}$ that fibers 
over it is
unique up to symplectomorphism independent of the choice of
$p_t$.
\end{prop}

For the proof of this we need Zung's classification of integrable Hamiltonian
systems with non-degenerate singularities, phrased in terms of
Lagrangian fibrations~\cite{Zung.II}:

\begin{defn}\label{rough.def}
Two singular Lagrangian fibrations $\rho_i\co (M_i,\omega_i)\rightarrow B_i$,
$i=1,2$, 
are {\it roughly symplectically equivalent} if there 
is an open cover $\{U_\alpha\}$ of $B_1$, a homeomorphism 
$\rho\co B_1\rightarrow B_2$, and
fiber preserving symplectomorphisms 
$\Phi_\alpha\co \rho_1^{-1}(U_\alpha)\rightarrow
\rho_2^{-1}(\phi(U_\alpha))$ such that on $\rho_1^{-1}(U_\alpha\cap U_\beta)$
the map $\Phi_\alpha^{-1}\circ\Phi_\beta$ induces the identity map on the 
fundamental group of the strata of each fiber and the identity map 
on the first integral homology of each fiber.
\end{defn}

Here the fibers are stratified as unions of orbits when one views
$\rho_1^{-1}(U_\alpha)$ as an integrable Hamiltonian system by composing
$\rho_1$ with a map $F\co U_\alpha\rightarrow R^n$.

\begin{thm}{\rm \cite{Zung.II}} \label{equiv.thm}
Two singular Lagrangian fibrations that are roughly symplectically equivalent
are fiberwise symplectomorphic if and only if 
they have the same Lagrangian class with respect to a common reference system.
\end{thm}

The Lagrangian class of $\pi\co (M,\omega)\rightarrow B$
is an element of $H^1(B,\cZ/\cR)$ where $\cZ$ is the
sheaf of local closed $1$-forms $\alpha$ on $M$ such that $\iota_X\alpha=
\iota_Xd\alpha$ for any vector $X$ such that $\pi_*X=0$ and $\cR$ is the
sheaf of symplectic fiber-preserving $S^1$ actions.
Identifying a reference Lagrangian fibration is necessary when there is no
roughly symplectically equivalent fibration that has a section.

\begin{proof}[Proof of Proposition~\ref{ball.prop}]
Cover the base $U_{n,m}$ with a collar neighborhood $V_b$
of the boundary and a disk neighborhood $V_{p_t}$ of $p_t$.
Then $V_b$ determines a unique Lagrangian fibered manifold~\cite{BouMol.fibr}
and by an isotopy such as in the proof of Lemma~\ref{nbhd.lem}
we can assume that $V_{p_t}$ determines a unique Lagrangian fibered manifold.
Choosing $\phi$ to be the identity map,
the conditions of Definition~\ref{rough.def} are met because the affine
structure on the base determines, up to isomorphism, 
the sublattice of $H_1(F,\bZ)$ generated by
the cycles of a regular fiber $F$ 
that collapse as the fiber moves to the boundary and to the nodal fiber.
Finally, because $H^1(U_{n,m},\cZ/\cR)=0$, Theorem~\ref{equiv.thm} implies
the fibrations are symplectically equivalent~\cite{Zung.II}.
\end{proof}

If we vary the position of $p_t$  (by varying our choice of $t$)
we get a family of symplectic
forms on the rational ball, all of which are equal near the boundary. 
Again, the vanishing of the rational cohomology
of $B_{n,m}$ allows
a Moser argument to confirm that the symplectic structures are isotopic.

The essential element for our  proof of Theorem~\ref{gbd.thm} is the fact that
a collar neighborhood of the boundary of $B_{n,m}$ 
is well defined up to fiberwise symplectomorphism by its base $V_b$.

\section{The symplectic surgery}

With the symplectic models for $B_{n,m}$ and $C_{n,m}$ at hand,
the proof of Theorem~\ref{gbd.thm} amounts to observing that 
we can choose $B_{n,m}$ and $C_{n,m}$ so that collar neighborhoods
of their boundaries symplectically 
embed into $L(n^2,nm-1)\times(0,\infty)$, fibering over $V_{n^2,nm-1}$
in such a way that their images in $V_{n^2,nm-1}$ coincide.

\smallskip
\begin{proof}[Proof of Theorem \ref{gbd.thm}] 
As explained at the end of Example~\ref{nbhd.ex},
given a symplectic $4$-manifold $(M,\omega)$ and an embedding 
$\psi\co C_{n,m}\rightarrow M$ 
such that each sphere $\psi(S_i)$ is a symplectic submanifold
we can assume the embedding $\psi$ is symplectic and gives
a Lagrangian fibration 
$\pi\co (\psi(C_{n,m}),\omega)\rightarrow W_{n,m}\subset\bR^2$.

Following Example~\ref{nbhd.ex} we can choose
$u_0=(0,-1)$ and $u_1=(1,0)$, so the vector $u_{k+1}$ defines a line in 
$\bR^2$ with slope $\frac{nm-1}{n^2}$.
Now $\psi(C_{n,m}-\cup_{i=1}^k S_k)$ fibers over 
$W_{n,m}'= W_{n,m}-\cup_{i=1}^k u_k$, so $W_{n,m}'$
defines a collar neighborhood of the boundary of $\psi(C_{n,m})$.
But $W_{n,m}'$ can also be viewed as a subset of $A_{n,m}$
so long as the distinguished point $(tn,tm)$ is chosen with
$t$ sufficiently small.  (See Section~\ref{rball.sec} and 
Figure~\ref{rball.fig}.)
As a subset of $A_{n,m}$ we see that $W_{n,m}'$ defines a collar neighborhood
of the boundary of a rational ball $B_{n,m}$.
Since these two collar neighborhoods fiber over the same simply connected
base they are symplectomorphic.  Therefore, we can find a symplectomorphism
$\phi$ that equips the generalized rational blowdown,
$\wtilde M=(M-\psi(\cup_{i=1}^k S_i))\cup_\phi  B_{n,m}$, 
with a symplectic
structure coming from those on $M$ and $B_{n,m}$.
\end{proof}

As for the rational blowdown with $m=1$, 
the volume of the generalized rational blowdown $\wtilde M$ 
is independent of any choice of rational
ball that fits.
The argument is exactly the same as in~\cite{Sym.blowdowns}.
It would be interesting to know whether a rational
blowdown, generalized or not, is unique up to symplectomorphism.

In the above proof we did not mention what is typically a crucial issue 
when trying to prove a surgery can be done symplectically: symplectic 
convexity of the neighborhood on which the gluing takes place.
A symplectic manifold $(M,\omega)$ with nonempty boundary
is {\it symplectically convex} if
there is an expanding vector field $X$ defined near and transverse to 
the boundary.  
To say that $X$ is expanding means $X$ points outward and 
$\cL_X \omega=\omega$.
The expanding vector field $X$ defines a contact structure on the
boundary, the $2$-plane field defined 
as the kernel of the $1$-form $\iota_X\omega$ restricted to the boundary.

If $A$ and $B$ are symplectic $2n$-manifolds with 
contactomorphic symplectically
convex boundaries and $A\subset (M,\omega)$ where $M$ is $2n$-dimensional,
then $(M-{\rm int}\ A)\cup  B $ admits
a symplectic structure induced from those of $M$ and $B$
(See~\cite{Etnyre.convexity} for more about symplectic convexity, contact
structures and symplectic surgeries.)

Thanks to the model spaces, we get symplectic convexity and contactomorphic
boundaries for free as follows.
Using the same notation as in the proof of Theorem~\ref{gbd.thm}, 
we can choose an arbitrarily small Lagrangian fibered neighborhood of
spheres $C_{n,m}$ fibering over a $W_{n,m}$ such that  the boundary of
$W'_{n,m}$ is transverse to the vector field 
$p_1 \frac{\del}{\del p_1}+p_2 \frac{\del}{\del p_2}$
(when viewed as a subset of $V_{n^2,nm-1}$).
This vector field lifts an expanding vector field on the preimage
of $W'_{n,m}$, thereby demonstrating the symplectic convexity of $C_{n,m}$.
Since we construct $B_{n,m}$ so that a collar neighborhood of its boundary
is symplectomorphic to that of $C_{n,m}$, the contact equivalence of the
boundaries and the symplectic convexity of $B_{n,m}$ are immediate.


\Addresses
\recd
\end{document}

%% file: agtout.tex
\def\ifplaintex{\expandafter\ifx\csname documentclass\endcsname\relax}

\def\gtp{{\mathsurround=0pt\it $\cal G\mskip-2mu$eometry \&\ 
$\cal T\!\!$opology $\cal P\!$ublications}}  

\def\recd{{\small Received:\qua\receiveddate\ifx\reviseddate\relax
\else\qquad Revised:\qua\reviseddate\fi\par}} 

\def\lognumber#1{\def\thelognumber{#1}}
\def\volumenumber#1{\def\thevolumenumber{#1}}
\def\volumeyear#1{\def\thevolumeyear{#1}}
\def\papernumber#1{\def\thepapernumber{#1}}
\def\pagenumbers#1#2{\def\startpage{#1}\def\finishpage{#2}}
\def\published#1{\def\publishdate{#1}}

\def\received#1{\def\receiveddate{#1}}
\def\revised#1{\def\reviseddate{#1}}
\def\accepted#1{\def\accepteddate{#1}}

\long\def\asciiabstract#1{\long\def\theasciiabstract{#1}}

\let\\\par\let\thelognumber\relax\let\thevolumenumber\relax
\let\thepapernumber\relax\let\thevolumeyear\relax\let\startpage\relax
\let\finishpage\relax\let\publishdate\relax\let\receiveddate\relax
\let\reviseddate\relax\let\accepteddate\relax\let\theasciititle\relax
\let\theasciiauthors\relax
\let\theasciiabstract\relax

\let\theasciiemail\relax

\ifplaintex
\font\logobig=cmssbx10 scaled 3836
\font\logomed=cmssbx10 scaled 2557
\else
\font\logobig=cmssbx10 scaled 4200
\font\logomed=cmssbx10 scaled 2800
\fi

\long\def\makeagttitle{   
\count0=\startpage
\agt\hfill      
\hbox to 45truept{\vbox to 0pt{\vglue -13truept{\logomed A\kern -.37em{\logobig 
T}\kern -.38em G}\vss}\hss}
\break
{\small Volume \thevolumenumber\ (\thevolumeyear)
\startpage--\finishpage\nl
Published: \publishdate}

\vglue .25truein

{\parskip=0pt\leftskip 0pt plus
1fil\def\\{\par\smallskip}{\Large\bf\thetitle}\par\medskip} \vglue
0.05truein

{\parskip=0pt\leftskip 0pt plus 1fil\def\\{\par}{\sc\theauthors}
\par\medskip}%
 
\vglue 0.03truein 

{\small\leftskip 25truept\rightskip 25truept{\bf Abstract}\stdspace\theabstract

{\bf AMS Classification}\stdspace\theprimaryclass
\ifx\thesecondaryclass\relax\else; \thesecondaryclass\fi\par
{\bf Keywords}\stdspace \thekeywords\par}\vglue 7truept

}   

\ifplaintex
\font\phead=cmsl9 scaled 950
\font\pnum=cmbx10 scaled 913
\font\pfoot=cmsl9 scaled 950
\headline{\vbox to 0pt{\vskip -4.5mm\line{\small\phead\ifnum
\count0=\startpage ISSN 1472-2739 (on-line) 1472-2747 (printed)
\hfill {\pnum\folio}\else\ifodd\count0\def\\{ }%
\ifx\theshorttitle\relax\thetitle\else\theshorttitle\fi\hfill{\pnum\folio}
\else\def\\{ and }{\pnum\folio}\hfill\ifx\theshortauthors\relax\theauthors
\else\theshortauthors\fi\fi\fi}\vss}}
\footline{\vbox to 0pt{\vglue 0mm\line{\small\pfoot\ifnum\count0=\startpage
\copyright\ \gtp\hfill\else
\agt, Volume \thevolumenumber\ (\thevolumeyear)\hfill\fi}\vss}}
\else
\font=cmsl9 scaled 1050
\font\lnum=cmbx10 
\font=cmsl9 scaled 1050
\makeatletter
\def\@oddhead{{\small\lhead\ifnum\count0=\startpage ISSN 1472-2739 
(on-line) 1472-2747 (printed)\hfill {\lnum\number\count0}\else\ifodd\count0
\def\\{ }\ifx\theshorttitle\relax \thetitle \else\theshorttitle\fi\hfill
{\lnum\number\count0}\else\def\\{ and }{\lnum\number\count0}
\hfill\ifx\theshortauthors\relax 
\theauthors\else\theshortauthors\fi\fi\fi}}\def\@evenhead{\@oddhead}
\def\@oddfoot{\small\lfoot\ifnum\count0=\startpage\copyright\ \gtp\hfill\else
\agt, Volume \thevolumenumber\ (\thevolumeyear)\hfill\fi}
\def\@evenfoot{\@oddfoot}
\makeatother
\fi
\let\maketitlepage\makeagttitle
\let\makeshorttitle\maketitlepage
\let\maketitle\maketitlepage

\newwrite\gtoutfile
\long\gdef\makeheadfile{  
{\def\\{, }\def\s{ }
\immediate\openout\gtoutfile head.xxx
\immediate\write\gtoutfile{To: math@arxiv.org}
\immediate\write\gtoutfile{Subject: put OR rep NNNNN:ppppp}
\immediate\write\gtoutfile{--text follows this line--}
\immediate\write\gtoutfile{Proxy-for: \ifx\theasciiauthors\relax
\theauthors\else\theasciiauthors\fi\s<\ifx\theasciiemail\relax\theemail\else\theasciiemail\fi>}
\immediate\write\gtoutfile{\noexpand\\}
\immediate\write\gtoutfile{Authors: \ifx\theasciiauthors\relax
\theauthors\else\theasciiauthors\fi}
{\def\\{ }\immediate\write\gtoutfile{Title: \ifx\theasciititle\relax
\thetitle\else\theasciititle\fi}}
\immediate\write\gtoutfile{Subj-class: GT or SG, GR etc}
\immediate\write\gtoutfile{MSC-class: \theprimaryclass\ifx\thesecondaryclass\relax\else, \thesecondaryclass\fi}
\immediate\write\gtoutfile{Journal-ref: Algebr. Geom. Topol. \thevolumenumber\s
(\thevolumeyear) \startpage-\finishpage}
\immediate\write\gtoutfile{Comments: Published by Algebraic and
Geometric Topology at}
\immediate\write\gtoutfile{\s\s\s  http://www.maths.warwick.ac.uk/agt/AGTVol\thevolumenumber/agt-\thevolumenumber-\thepapernumber.abs.html}
\immediate\write\gtoutfile{\noexpand\\}
\immediate\write\gtoutfile{}
\ifx\theasciiabstract\relax
\immediate\write\gtoutfile{\theabstract}\else
\immediate\write\gtoutfile{\theasciiabstract}\fi
\immediate\write\gtoutfile{}
\immediate\write\gtoutfile{\noexpand\\}
\immediate\write\gtoutfile{}
\immediate\closeout\gtoutfile}}  

\def\maketitlepage{\makeagttitle\makeheadfile}
\let\makeshorttitle\maketitlepage
\let\maketitle\maketitlepage

\def\ifplaintex{\expandafter\ifx\csname documentclass\endcsname\relax}

\def\gtp{{\mathsurround=0pt\it $\cal G\mskip-2mu$eometry \&\ 
$\cal T\!\!$opology $\cal P\!$ublications}}  

\def\recd{{\small Received:\qua\receiveddate\ifx\reviseddate\relax
\else\qquad Revised:\qua\reviseddate\fi\par}} 

\def\lognumber#1{\def\thelognumber{#1}}
\def\volumenumber#1{\def\thevolumenumber{#1}}
\def\volumeyear#1{\def\thevolumeyear{#1}}
\def\papernumber#1{\def\thepapernumber{#1}}
\def\pagenumbers#1#2{\def\startpage{#1}\def\finishpage{#2}}
\def\published#1{\def\publishdate{#1}}

\def\received#1{\def\receiveddate{#1}}
\def\revised#1{\def\reviseddate{#1}}
\def\accepted#1{\def\accepteddate{#1}}

\long\def\asciiabstract#1{\long\def\theasciiabstract{#1}}

\let\\\par\let\thelognumber\relax\let\thevolumenumber\relax
\let\thepapernumber\relax\let\thevolumeyear\relax\let\startpage\relax
\let\finishpage\relax\let\publishdate\relax\let\receiveddate\relax
\let\reviseddate\relax\let\accepteddate\relax\let\theasciititle\relax
\let\theasciiauthors\relax
\let\theasciiabstract\relax

\let\theasciiemail\relax

\ifplaintex
\font\logobig=cmssbx10 scaled 3836
\font\logomed=cmssbx10 scaled 2557
\else
\font\logobig=cmssbx10 scaled 4200
\font\logomed=cmssbx10 scaled 2800
\fi

\long\def\makeagttitle{   
\count0=\startpage
\agt\hfill      
\hbox to 45truept{\vbox to 0pt{\vglue -13truept{\logomed A\kern -.37em{\logobig 
T}\kern -.38em G}\vss}\hss}
\break
{\small Volume \thevolumenumber\ (\thevolumeyear)
\startpage--\finishpage\nl
Published: \publishdate}

\vglue .25truein

{\parskip=0pt\leftskip 0pt plus
1fil\def\\{\par\smallskip}{\Large\bf\thetitle}\par\medskip} \vglue
0.05truein

{\parskip=0pt\leftskip 0pt plus 1fil\def\\{\par}{\sc\theauthors}
\par\medskip}%
 
\vglue 0.03truein 

{\small\leftskip 25truept\rightskip 25truept{\bf Abstract}\stdspace\theabstract

{\bf AMS Classification}\stdspace\theprimaryclass
\ifx\thesecondaryclass\relax\else; \thesecondaryclass\fi\par
{\bf Keywords}\stdspace \thekeywords\par}\vglue 7truept

}   

\ifplaintex
\font\phead=cmsl9 scaled 950
\font\pnum=cmbx10 scaled 913
\font\pfoot=cmsl9 scaled 950
\headline{\vbox to 0pt{\vskip -4.5mm\line{\small\phead\ifnum
\count0=\startpage ISSN 1472-2739 (on-line) 1472-2747 (printed)
\hfill {\pnum\folio}\else\ifodd\count0\def\\{ }%
\ifx\theshorttitle\relax\thetitle\else\theshorttitle\fi\hfill{\pnum\folio}
\else\def\\{ and }{\pnum\folio}\hfill\ifx\theshortauthors\relax\theauthors
\else\theshortauthors\fi\fi\fi}\vss}}
\footline{\vbox to 0pt{\vglue 0mm\line{\small\pfoot\ifnum\count0=\startpage
\copyright\ \gtp\hfill\else
\agt, Volume \thevolumenumber\ (\thevolumeyear)\hfill\fi}\vss}}
\else
\font=cmsl9 scaled 1050
\font\lnum=cmbx10 
\font=cmsl9 scaled 1050
\makeatletter
\def\@oddhead{{\small\lhead\ifnum\count0=\startpage ISSN 1472-2739 
(on-line) 1472-2747 (printed)\hfill {\lnum\number\count0}\else\ifodd\count0
\def\\{ }\ifx\theshorttitle\relax \thetitle \else\theshorttitle\fi\hfill
{\lnum\number\count0}\else\def\\{ and }{\lnum\number\count0}
\hfill\ifx\theshortauthors\relax 
\theauthors\else\theshortauthors\fi\fi\fi}}\def\@evenhead{\@oddhead}
\def\@oddfoot{\small\lfoot\ifnum\count0=\startpage\copyright\ \gtp\hfill\else
\agt, Volume \thevolumenumber\ (\thevolumeyear)\hfill\fi}
\def\@evenfoot{\@oddfoot}
\makeatother
\fi
\let\maketitlepage\makeagttitle
\let\makeshorttitle\maketitlepage
\let\maketitle\maketitlepage

\newwrite\gtoutfile
\long\gdef\makeheadfile{  
{\def\\{, }\def\s{ }
\immediate\openout\gtoutfile head.xxx
\immediate\write\gtoutfile{To: math@arxiv.org}
\immediate\write\gtoutfile{Subject: put OR rep NNNNN:ppppp}
\immediate\write\gtoutfile{--text follows this line--}
\immediate\write\gtoutfile{Proxy-for: \ifx\theasciiauthors\relax
\theauthors\else\theasciiauthors\fi\s<\ifx\theasciiemail\relax\theemail\else\theasciiemail\fi>}
\immediate\write\gtoutfile{\noexpand\\}
\immediate\write\gtoutfile{Authors: \ifx\theasciiauthors\relax
\theauthors\else\theasciiauthors\fi}
{\def\\{ }\immediate\write\gtoutfile{Title: \ifx\theasciititle\relax
\thetitle\else\theasciititle\fi}}
\immediate\write\gtoutfile{Subj-class: GT or SG, GR etc}
\immediate\write\gtoutfile{MSC-class: \theprimaryclass\ifx\thesecondaryclass\relax\else, \thesecondaryclass\fi}
\immediate\write\gtoutfile{Journal-ref: Algebr. Geom. Topol. \thevolumenumber\s
(\thevolumeyear) \startpage-\finishpage}
\immediate\write\gtoutfile{Comments: Published by Algebraic and
Geometric Topology at}
\immediate\write\gtoutfile{\s\s\s  http://www.maths.warwick.ac.uk/agt/AGTVol\thevolumenumber/agt-\thevolumenumber-\thepapernumber.abs.html}
\immediate\write\gtoutfile{\noexpand\\}
\immediate\write\gtoutfile{}
\ifx\theasciiabstract\relax
\immediate\write\gtoutfile{\theabstract}\else
\immediate\write\gtoutfile{\theasciiabstract}\fi
\immediate\write\gtoutfile{}
\immediate\write\gtoutfile{\noexpand\\}
\immediate\write\gtoutfile{}
\immediate\closeout\gtoutfile}}  

\def\maketitlepage{\makeagttitle\makeheadfile}
\let\makeshorttitle\maketitlepage
\let\maketitle\maketitlepage